\documentclass[12pt,reqno]{amsart}

\usepackage[cp1251]{inputenc}

\usepackage{amssymb
,epsfig
}
\usepackage{graphicx}

\newtheorem{thm}{Theorem}
\theoremstyle{definition}
\newtheorem{prp}{Proposition}
\newtheorem{lemma}{Lemma}

\newtheorem{ex}{Example}
\newtheorem{rem}{Remark}
\newtheorem{cor}{Corollary}

\newtheorem{dfn}{Definition}
\newtheorem{alg}{Algorithm}
\begin{document}

\title{Cyclopermutohedron
\newline
(To Nilolai Dolbilin, who loves polytopes)}

\author{ Gaiane Panina}

\address{*
 $^\dag$ Institute for Informatics and Automation, St. Petersburg, Russia,
Saint-Petersburg State University, St. Petersburg, Russia,
e-mail:gaiane-panina@rambler.ru.  }


\begin{abstract}
It is known that the $k$-faces of the permutohedron $\Pi_n$ can be
labeled by (all possible)  linearly ordered partitions of the set
$[n]=\{1,...,n\}$  into $(n-k)$  non-empty parts. The incidence
relation corresponds to the refinement: a face $F$ contains a face
$F'$ whenever the  label of $F'$ refines
 the  label of $F$.

 In the paper we consider the cell complex  ${CP}_{n+1}$ defined in analogous way, replacing linear ordering by cyclic ordering.
  Namely, $k$-cells of the complex ${CP}_{n+1}$ are labeled by (all possible)  cyclically ordered
partitions of the set $[n+1]=\{1,...,n, n+1\}$  into $(n+1-k)$
non-empty parts, where $(n+1-k)>2$. The incidence relation in
${CP}_{n+1}$ again corresponds to the refinement: a cell $F$
contains a cell $F'$ whenever the  label of $F'$ refines
 the  label of $F$.

 The complex ${CP}_{n+1}$ cannot be represented by a convex polytope, since   it is not a combinatorial sphere (not even a combinatorial manifold).
 However, it can be represented by some  \textit{virtual polytope} (that is, Minkowski difference of two convex polytopes)
 which we call "cyclopermutohedron"  $\mathcal{CP}_{n+1}$. It is defined explicitly, as a weighted Minkowski sum of line segments.
 Informally, the cyclopermutohedron can be viewed as "permutohedron with diagonals".
One of the motivations is that the cyclopermutohedron is  a
"universal" polytope for moduli spaces of polygonal linkages.
\end{abstract}

\maketitle

\section{Preliminaries}\label{SectionPrelim}
\begin{figure}[h]\label{permutahedron}
\centering
\includegraphics[width=6 cm]{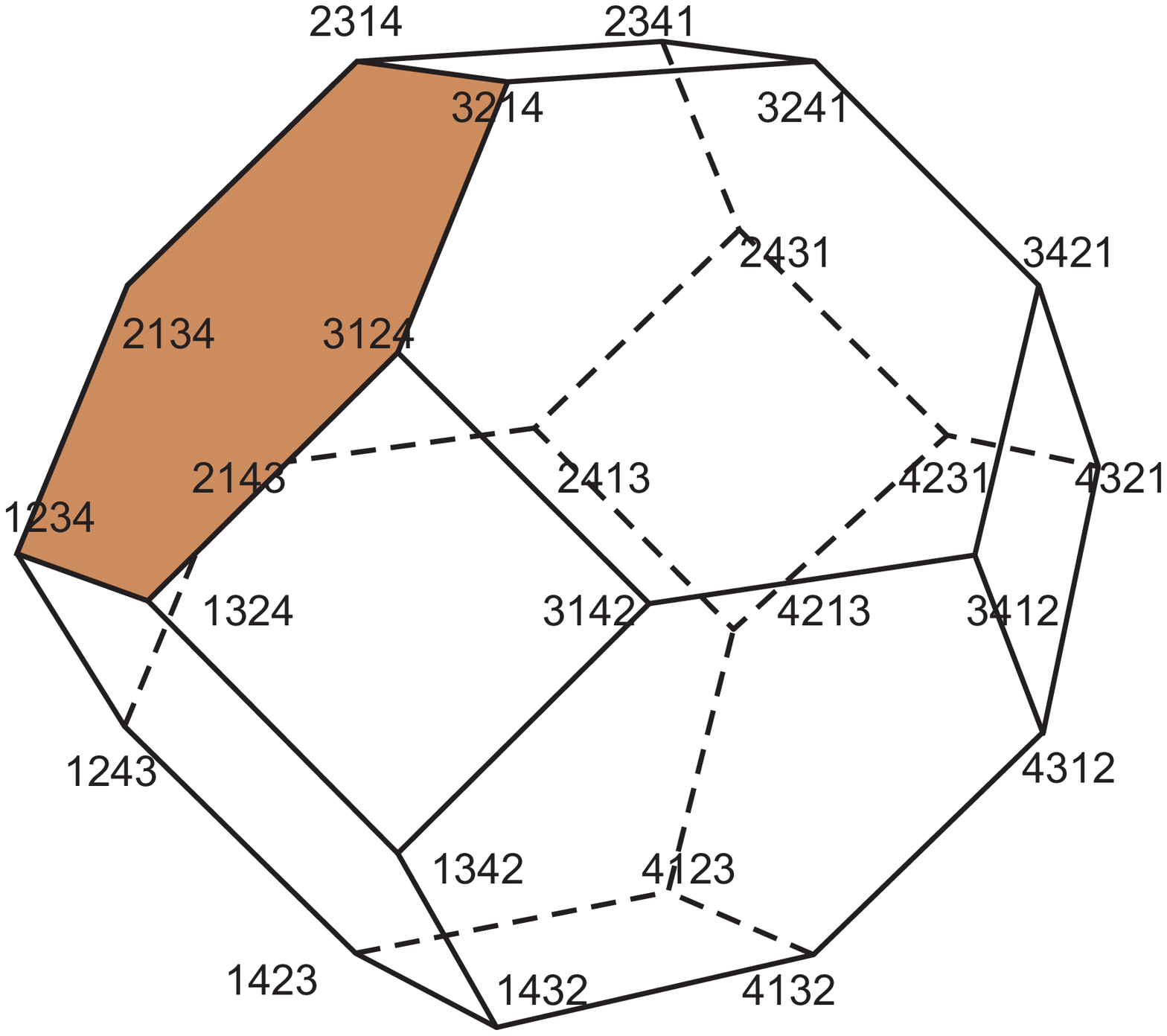}
\caption{Permutahedron $\Pi_4$ with labels on the vertices. The
shadowed face is labeled by $(\{1,2,3\}\{4\})$.}
\end{figure}

Among all famous polytopes with  combinatorial backgrounds, the
\textit{standard permutohedron} $\Pi_n$ is the oldest and probably
the most important one. It  is defined (see \cite{z}) as the convex
hull of all points in $\mathbb{R}^n$ that are obtained by permuting
the coordinates of the point $(1,2,...,n)$. It has the following
properties:
\begin{enumerate}
    \item $\Pi_n$ is an $(n-1)$-dimensional polytope.
    \item The $k$-faces of  $\Pi_n$
are labeled by ordered partitions of the set
\newline $[n]=\{1,2,...,n\}$  into
$(n-k)$ non-empty parts.

\item In particular, the vertices are labeled by the elements of the
symmetry group $S_n$. For each vertex, the label is the inverse
permutation of the coordinates of the vertex. Two vertices are
joined by an edge whenever their labels differ by a permutation of
two neighbor entries.

    \item A face $F$ of $\Pi_n$ is contained in a face $F'$  iff the label of
$F$ is a refinement of the label  of $F'$.

Here and in the sequel, we mean the order-preserving refinement. For
instance, $(\{1,3\},\{5,6\},\{4\},\{2\})$ refines
$(\{1,3\},\{5,6\},\{2,4\})$, but does not refine
$(\{1,3\},\{2,4\},\{5,6\})$.
    \item Each face of $\Pi_n$ equals the Cartesian product of standard permutohedra of smaller
dimensions.
    \item The permutohedron is a \textit{zonotope}, that is,  Minkowski sum
of line segments.
\item
 The permutohedra $\Pi_1$, $\Pi_2$, and $\Pi_3$ are a one-point polytope,  a
 segment,
 and a regular hexagon respectively.  The permutohedron $\Pi_4$ is depicted in Figure \ref{permutahedron}.
\end{enumerate}

Historically, the permutohedron was followed by associahedron,
cyclohedron  (see \cite{z}), permutoassociahedron  (see \cite{K}),
and more advanced polytopes  appearing first as combinatorially
described cell complexes: secondary polytope \cite{GKZ}, generalized
associahedra \cite{Fomin}. The complexes luckily appeared to be
combinatorial spheres, and moreover, representable by convex
polytopes. A wider framework includes graph-associahedra, nestohedra
\cite{Post}, and other recently constructed polytopes.

In the paper, we proceed in the same manner: we first describe a
cell complex ${CP}_{n+1}$, and then show that it can be represented
by a \textit{virtual polytope},
 i.e., by  Minkowski difference of two convex polytopes. We guess that virtual polytopes  (as the tightest algebraic generalization of convex
  polytopes, originally introduced in \cite{pkh}) provide us a reasonable framework for this particular complex.

The combinatorics of the complex ${CP}_{n+1}$ almost literally
repeats the combinatorics of permutohedron with one essential
difference: linear ordering is replaced by cyclic ordering.
 Thus defined, complex ${CP}_{n+1}$ cannot be represented by a convex polytope, since it is not a combinatorial sphere (not even a combinatorial manifold). However, it can be represented by a virtual polytope
 which we call \textit{cyclopermutohedron}  $\mathcal{CP}_{n+1}$. It is defined explicitly, as a weighted Minkowski sum of line segments.
 The word "weighted" means that  some of the summands are negatively weighted, that is, we take not only Minkowski sum, but also\textit{ Minkowski  difference}.

In an oversimplified way, the cyclopermutohedron
$\mathcal{CP}_{n+1}$  can be visualized as the permutohedron $\Pi_n$
"with diagonals". This means that all the  proper faces of  $\Pi_n$
are also faces of  $\mathcal{CP}_{n+1}$. However,
$\mathcal{CP}_{n+1}$ has some extra faces in comparison with
$\Pi_n$. The latter look like "diagonal" faces, see Figure
\ref{FigCycloper5labels}.

The idea to replace linear ordering by cyclic ordering in the
framework of "famous polytopes"  is not new:   replacing linear
ordering in the combinatorics of  associahedron  yields the
cyclohedron.

Although the cyclopermutohedron seems to be interesting  for its own
sake, we were initially motivated by yet another application: the
cyclopermutohedron  is a "universal" polytope for moduli spaces of
polygonal linkages. That is, for any length assignment of a planar
flexible polygon  with $(n+1)$ edges, its moduli space admits a
natural cell structure which embeds in the (cell structure of) the
cyclopermutohedron $\mathcal{CP}_{n+1}$.

The paper is organized as follows. Section \ref{SecCombComplex}
describes a regular cell complex which we wish to realize as the
face lattice of the cyclopermutohedron. Section
\ref{SecCyclopermVIP} contains the construction, the main statements
(Proposition \ref{ProposVert}  and Theorem 1), and low dimensional
examples of the cyclopermutohedra. The main theorem (Theorem 1)
establishes combinatorial equivalence of the cyclopermutohedron and
the cell complex.  Next, Section \ref{SecModuliPolygon} explains the
relationship with the moduli space of a planar polygonal linkage.
The proofs of Theorem 1  and Proposition \ref{ProposVert} are in
Section \ref{SectProof}. To make the paper self-contained, we put
all necessary information on virtual polytopes   in Section
\ref{SecVIPs}.

\section{Cyclopermutohedron: the desired combinatorics}\label{SecCombComplex}

By a \textit{regular cell complex } we mean a finite cell complex
which can be constructed inductively by defining its skeleta. The
$0$-skeleton is a finite set of points endowed with discrete
topology. Once the $(k - 1)$-skeleton is constructed, we attach a
collection of closed $k$-dimensional balls $C_i$, called
\textit{closed cells}. For each of $C_i$, the attaching map is a
homeomorphism between $\partial C_i$ and a sphere homeomorphic
subcomplex of the $(k - 1)$-skeleton.

By construction, the boundary of each closed cell has the induced
structure of a regular cell complex.

To define a regular cell complex, it suffices to list all the closed
cells of the complex together with the incidence relations.
Following this rule, for a fixed   number $n>2$, we define the
regular cell complex ${CP}_{n+1}$
 as follows.
 \begin{itemize}
   \item For  $k=0,...,n-2$, the   $k$-dimensional cells ($k$-cells, for short) of the complex ${CP}_{n+1}$ are labeled by (all possible)  cyclically ordered
partitions of the set $[n+1]=\{1,...,n+1\}$  into $(n-k+1)$  non-empty
parts.
   \item A (closed) cell $F$ contains a cell $F'$ whenever the  label of $F'$ refines
 the  label of $F$. Here  we mean orientation
 preserving refinement.
 \end{itemize}

In particular, this means that the vertices of the complex
${CP}_{n+1}$ are labeled by cyclic orderings on the set $[n+1]$. Two
vertices are joined by an edge whenever their labels differ on a
permutation of two neighbor entries.

\bigskip

Let us first give an example, next agree how we depict the labels,
and then prove the correctness of the above construction.

\begin{figure}[h]
\centering
\includegraphics[width=8 cm]{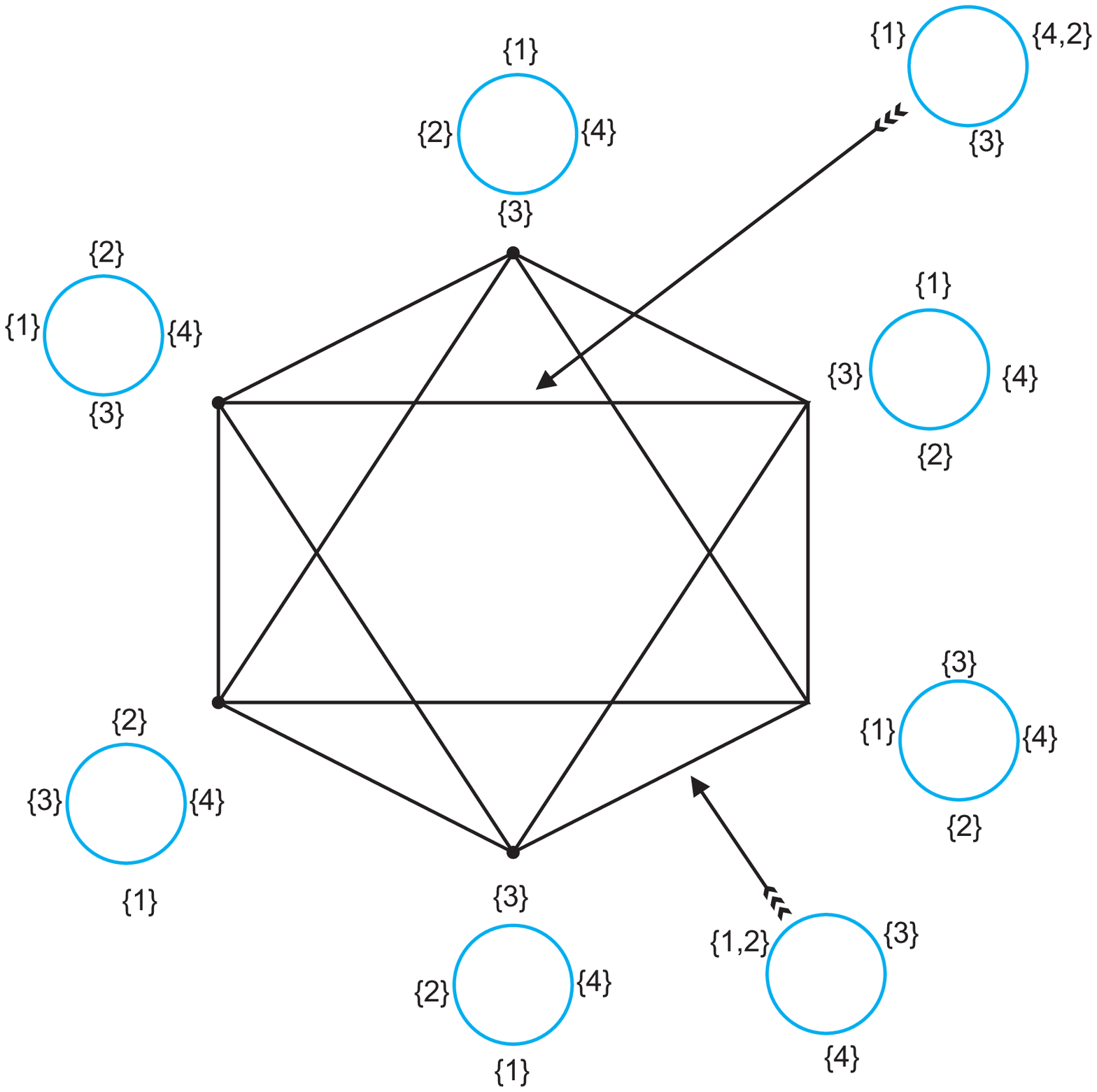}
\caption{The complex ${CP}_{4}$. }\label{Fig1}
\end{figure}

\begin{ex} For $n=3$, the complex  ${CP}_{n+1}={CP}_{4}$ consists of  six vertices
 and twelve edges. In the figure \ref{Fig1}, all the vertices and
two of the edges are labeled.
\end{ex}

\subsection*{A convention about presenting the labels}
Let us agree that instead of saying "the cell of the complex labeled
by $\lambda$" we say for short "the cell  $\lambda$".

We shall write a cyclically ordered partition $\lambda$ in the
following ways:
\begin{enumerate}
\item Sometimes  we depict $\lambda$ as a counterclockwise oriented circle with subsets placed on the circle, as in Figure \ref{Fig1}.
    \item Alternatively, we cut  $\lambda$ at any place and  get a string of sets $(I_1,I_2,...,I_k)$.
    We keep in mind that $$(I_1,I_2,...,I_k)=cycl(I_1,I_2,...,I_k)=(I_k,I_1, I_2,...,I_{k-1})=$$
    $$=cycl^2(I_1,I_2,...,I_k)=(I_{k-1},I_k, I_1,I_2, ...,I_{k-3},I_{k-2}), \hbox{ etc}.$$
    \item For the vertices of  ${CP}_{n+1}$, we use one more type of reduction. Given a label  $\lambda$  of a vertex
(that is, a cyclic ordering on $[n+1]$), we write its label as the
linear ordering on $[n]$, which is obtained from $\lambda$  by
cutting through the entry $(n+1)$ and removing the entry  $(n+1)$
from the string. The obtained label we denote by
$\overline{\lambda}$.

For instance, the label
$$\lambda=(\{3\},\{4\},\{1\},\{5\},\{6\},\{2\})$$ yields
$$\overline{\lambda}=(2,3,4,1,5).$$

Conversely, given a linear ordering $\overline{\lambda}$ on the set
$[n]$, we cook a cyclically
 ordered partitions of the set $[n+1]$ into $(n+1)$ singletons  by adding at the end the set $\{n+1\}$ and closing the string in a cycle.

 For instance, $(3,2,1)\in S_3$ yields $(\{3\},\{2\},\{1\}),$ and then the cyclically ordered partition
 $$  \begin{array}{ccc}
        & \{3\} &  \\
       \{2\} &  & \{4\}. \\
        & \{1\} &
     \end{array} $$
\end{enumerate}

\begin{prp} For each $n\geq 3$, the  complex ${CP}_{n+1}$ is defined correctly and uniquely.
\end{prp}
Proof. We construct the complex inductively. The $0$-skeleton is
obviously well-defined. Assume that we have already the $k$-skeleton
and that we wish to attach a $(k+1)$-cell $C$ labeled by some
cyclically ordered partition
$$(I_1,I_2,...,I_{n-k})=(I_2,I_3,...,I_{n-k},I_1)= (I_3,I_4,...,I_1,I_2)=....$$ All the cells of
the $k$-skeleton that should be incident to $C$ form a subcomplex of
the $k$-skeleton which is combinatorially isomorphic to the boundary
complex of the Cartesian product of permutohedra
$$\partial(\Pi_{n_1}\times...\times\Pi_{n_{n-k}}),$$
where $n_i=|I_i|$. Obviously, it is a  $k$-dimensional sphere, and
there is a unique way to   attach the boundary $\partial C$ to the
sphere.\qed

\begin{cor}
Each closed cell  of the complex  ${CP}_{n+1}$ is combinatorially
equivalent to a Cartesian product of permutohedra. \qed
\end{cor}

\begin{lemma}Let $\lambda=(I_1,I_2,...,I_k)$ be the label of some cell  of
the  complex $ {CP}_{n+1}$. We assume that the entry $(n+1)$ belongs to the set $I_k$.

The (labels of the) vertices of the cell $\lambda$ can be retrieved
via the following algorithm:
\begin{enumerate}
    \item Take the string $\lambda=(I_1,I_2,...,I_k)$ and remove the entry   $n+1$ from $I_k$.
    \item Take all possible orderings inside each of the sets $I_j, \ j=1,2,...,k,$ and
    list all the resulting strings. We get some  set $V(\lambda)$.
    \item Apply  cyclic permutations   $cycl^0,cycl^1,  cycl^2,...,
    cycl^{|I_k|-1}$ to all the elements of $V(\lambda)$. We get
    $$V_0(\lambda)=V(\lambda), \ V_1(\lambda):=cycl^1V(\lambda),...,V_{|I_k|-1}(\lambda):=cycl^{|I_k|-1}V(\lambda).$$
\item The vertex set of the cell $\lambda$ is the union of the sets: $$ Vert(\lambda)= V_0\bigcup
V_1\bigcup V_2\bigcup...\bigcup V_{|I_k|-1}.\qed$$
\end{enumerate}
\end{lemma}

\begin{ex}\label{ExVertices}
For $\lambda=(\{1,2\},\{3\},\{4,5,6\})$, we have

 $$ Vert(\lambda)= $$
$$\{ (1,2,3,4,5),(1,2,3,5,4), \ (2,1,3,5,4),\ (2,1,3,4,5),$$
$$(5,1,2,3,4), \ (4,1,2,3,5), \ (4,2,1,3,5),\ (5,2,1,3,4),$$
$$(4,5,1,2,3,), \ (5,4,1,2,3), \ (5,4,2,1,3),\ (4,5,2,1,3,) \}.\qed$$

\end{ex}

To give a reader more intuition, we complete the section by yet
another example.

\begin{figure}[h]
\centering
\includegraphics[width=8 cm]{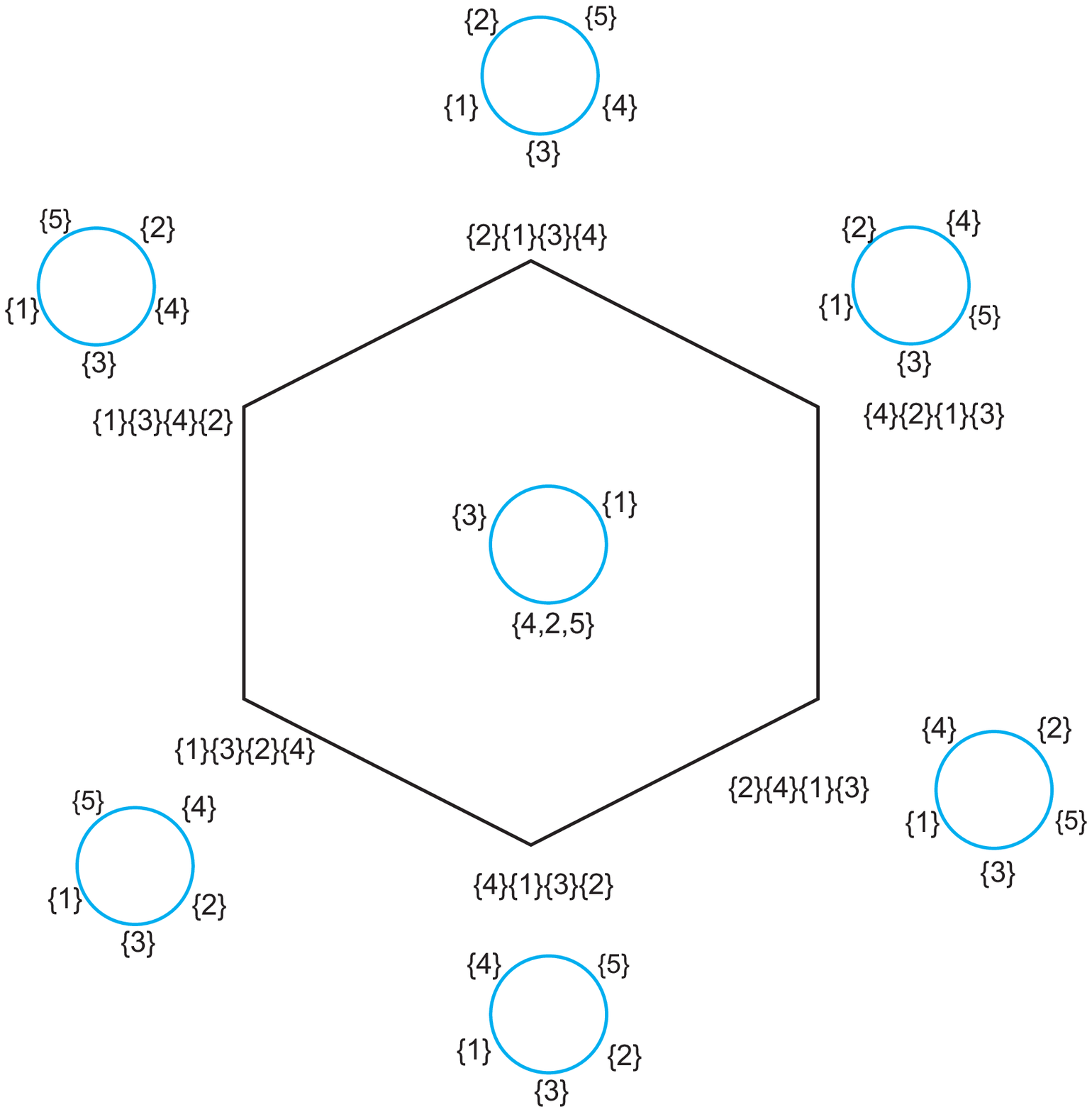}
\caption{The $2$-face  of the complex $CP_5$  labeled by
$(\{1\}\{3\}\{4,2,5\})$ }\label{Fig2}
\end{figure}

\begin{ex} For $n=4$,  Figure \ref{Fig2} presents  a hexagonal $2$-face  of the
complex $CP_5$  and all the vertices of the face.
\end{ex}

\section{Cyclopermutohedron:  realization as a virtual polytope}\label{SecCyclopermVIP}

Assuming that $\{e_i\}$  are standard basic vectors in $\mathbb{R}^n$, define the points

$$\begin{array}{ccccccccc}
    R_i=\sum_{j=1}^n (e_j-e_i)=(-1, & ... & -1, & n-1, & -1, & ... & -1, & -1, &-1, )\in \mathbb{R}^{n},\\
     &  & & \ i &  &  &  &  &
  \end{array}
$$
and  the following two families of line segments:
$$q_{ij}=\left[e_i,e_j\right], \ \ \  i<j$$ and $$  r_i=\left[0,R_{i} \right].$$

We also need the point $S=\left(1,1,...,1\right)\in \mathbb{R}^{n}$.

\begin{dfn} The  \textit{cyclopermutohedron} is a virtual polytope defined as    the weighted Minkowski sum
of line segments:
$$ \mathcal{CP}_{n+1}:=   \sum_{i< j} q_{ij} + S- \sum_{i=1}^n  r_i.$$

\end{dfn}
 Clearly, $ \mathcal{CP}_{n+1}$
 lies in the hyperplane
 $$x_1+...+x_n=\frac{n(n+1)}{2},$$
 so its actual dimension is $(n-1)$.

\begin{rem}\label{RemPermSum} The Minkowski sum
$$   \sum_{i< j} q_{ij}+ S$$  is known to be equal to the standard permutohedron
$\Pi_n$  (see \cite{z}). Therefore we can write

$$ \mathcal{CP}_{n+1}=  \Pi_n - \sum_{i=1}^n  r_i.$$

\end{rem}
\begin{rem}\label{symm}
The symmetric group $S_n$ acts on the $\mathbb{R}^n$ by permuting
the coordinates. The action preserves the sets $\{q_{i,j}\}$ and
$\{r_i\}$, so both permutohedron and cyclopermutohedron are
invariant under the action. We shall refer to this property as the
\textit{symmerty property} of the polytopes.
\end{rem}

A curios and crucial feature of  the cyclopermutohedron is the
following:

\begin{prp}\label{ProposVert} The set of vertices of $\mathcal{CP}_{n+1}$ equals the set of vertices of the standard permutohedron
$\Pi_n$:
$$Vert(\mathcal{CP}_{n+1})=Vert(\Pi_n).\qed$$

\end{prp}

\bigskip
The proposition allows us to label the vertices of the
cyclopermutohedron by the labels borrowed from the permutohedron.

More detailed, each vertex $v$ of $\mathcal{CP}_{n+1}$ is a vertex of $\Pi_n$,  and therefore,
 its coordinates is some permutation on $[n]$. We
label the vertex $v$ by the inverse permutation. For instance, the
vertex with coordinates  $(4,1,5,3,2)$ is labeled by $\lambda=(2,5,4,1,3)$.

Next, remind that each label yields automatically a cyclic ordering on the set $[n+1]$
by adding  to $\lambda$ the entry  $(n+1)$ at the end and closing the label in the circle.

Therefore we have:
\begin{cor}\label{CorNumberingVertices}The vertices  of
$\mathcal{CP}_{n+1}$ are canonically labeled by (all possible)
cyclic orderings on the set $[n+1]$.
\end{cor}

\begin{dfn} We  define a bijection
 $$\varphi: Vert(CP_{n+1})\rightarrow
Vert(\mathcal{CP}_{n+1})$$ between the vertex sets of the cell
complex $CP_n$ and the cyclopermutohedron. Given a vertex $v \in
Vert(CP_{n+1})$ labeled by $\lambda $, the bijection $\varphi$  maps
it to the vertex of the $ \mathcal{CP}_{n}$  labeled by
$\overline{\lambda }$. By construction, it coincides with the vertex
of $\Pi_n$ labeled by $\overline{\lambda }$.
\end{dfn}

\begin{thm}\label{TheoremMain}The   bijection $\varphi$  extends to a combinatorial
isomorphism between the complex $CP_{n+1}$  and the
cyclopermutohedron $ \mathcal{CP}_{n+1}$.

More precisely, the following holds:\begin{enumerate}
    \item  A subset  $ V \subset
Vert(CP_{n+1})$ is the vertex set of some cell of the complex
$CP_{n+1}$
 if and only if
 $\varphi( V)$ is the vertex set of some face of
$\mathcal{CP}_{n+1}$.
    \item The $k$-faces  of $\mathcal{CP}_{n+1}$ can be labeled by (all possible) cyclically ordered partitions of the set $[n+1]$  into $(n-k+1)$ non-empty parts.
    \item With this labeling, a face $F'$ is a face of a face $F$ if and only if the  label of $F'$ refines
 the  label of $F$.

\end{enumerate}
 \qed
\end{thm}

\begin{figure}[h]
\centering
\includegraphics[width=12 cm]{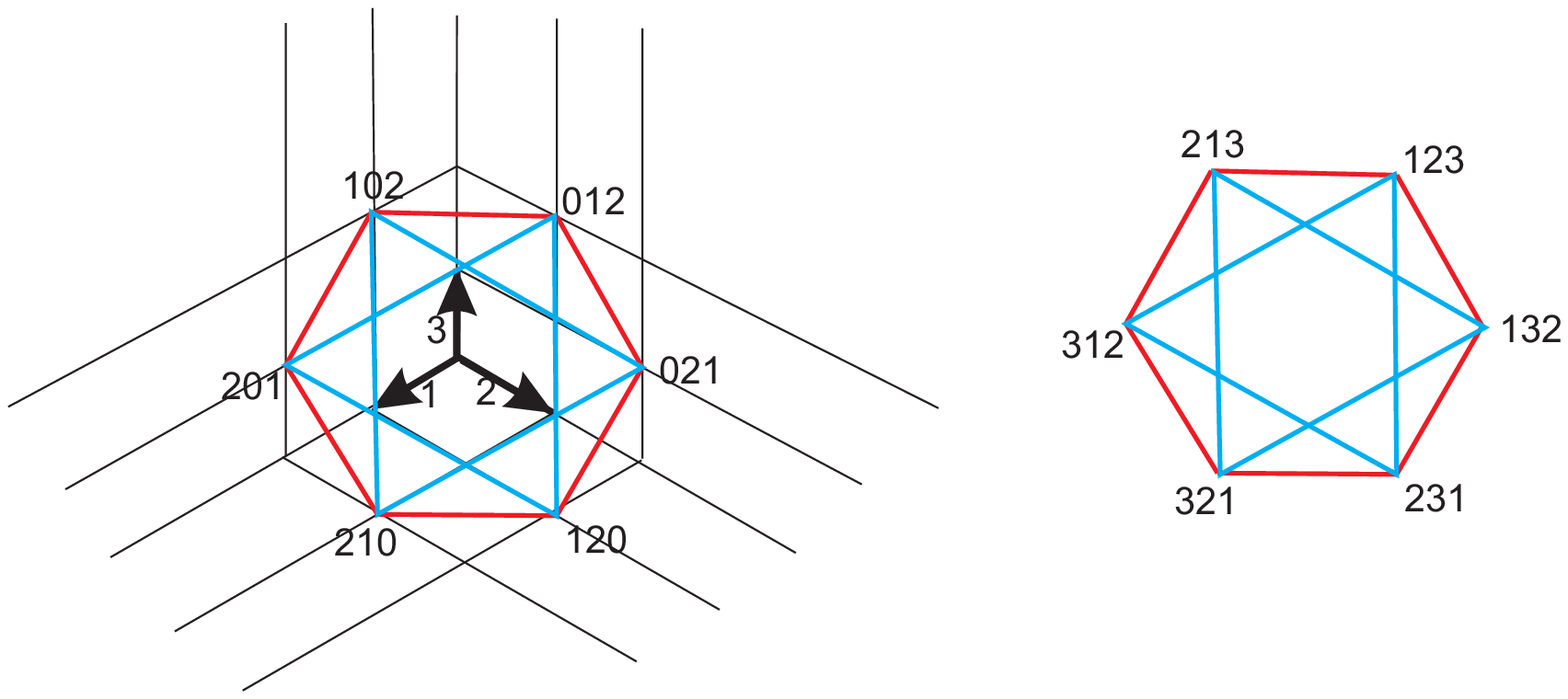}
\caption{Basis vectors in $\mathbb{R}^3$, the polytope
$\sum_{i< j} q_{ij} - \sum_{i=1}^3  r_i$
(left), and the cyclopermutohedron $\mathcal{CP}_{4}$ with
coordinates of the vertices (right).}\label{FigCycloper4}
\end{figure}

\begin{ex}\label{ExCycloper4} $ \mathcal{CP}_{4}$ is a two-dimensional virtual polytope, however, the definition puts in in $\mathbb{R}^3$.
 Figure \ref{FigCycloper4}, left  depicts the  virtual polytope
 $$  \sum_{1 \geq i> j \geq 3} q_{ij} - \sum_{i=1}^3  r_i  .$$ After the parallel translation by  $S=(1,1,1)$,
we get  the cyclopermutohedron $ \mathcal{CP}_{4}$ (right). (A more
detailed explanation of how to compute Minkowski difference in the
plane one finds in Section \ref{SecVIPs}). In the figure, we
indicate the coordinates of the vertices. After replacing
coordinates by labels, we arrive exactly at the complex $CP_4$
(compare with Figure \ref{Fig1}).
\end{ex}

\begin{figure}[h]
\centering
\includegraphics[width=10 cm]{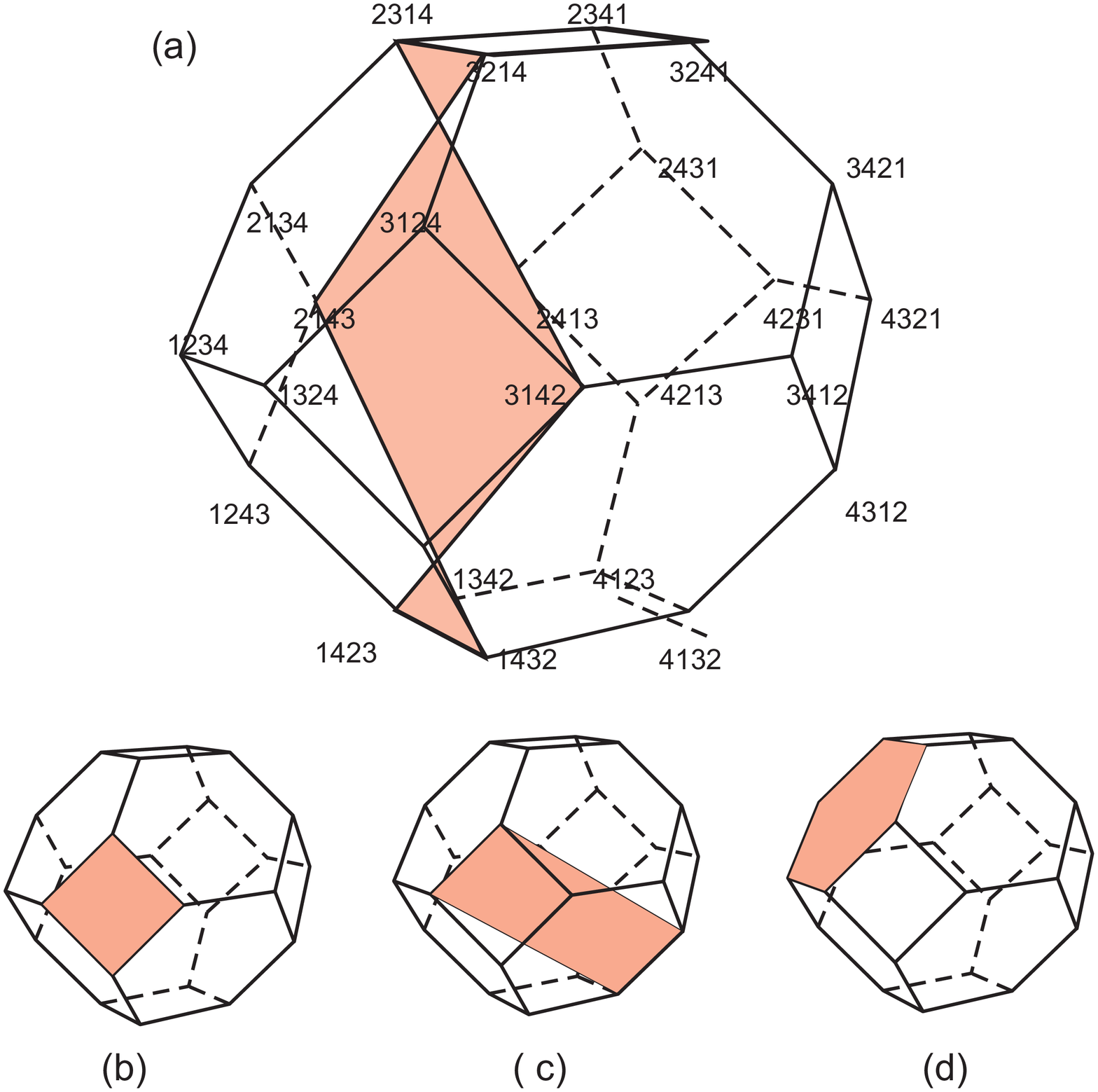}
\caption{There are four types of $2$-faces of the cyclopermutohedron
$ \mathcal{CP}_{5}$. We depict here representatives of all the
types. Their labels are: (a) $(\{1\}\{4\}\{2,3,5\})$, (b)
$(\{1,3\}\{4,2\}\{5\})$, \newline (c) $(\{1,3\}\{2\}\{4,5\})$, (d)
$(\{1,2,3\}\{4\}\{5\})$.}\label{FigCycloper5labels}
\end{figure}

\begin{ex} The
cyclopermutogedron $ \mathcal{CP}_{5}$  is a $3$-dimensional virtual
polytope. Here is the  complete list of its faces (see also Figure
\ref{FigCycloper5labels}).
\begin{enumerate}
  \item By Proposition \ref{ProposVert}, the  vertices of $ \mathcal{CP}_{5}$ are  the vertices  of the standard permutohedron $\Pi_4$.
  \item The edges are of two types: those coinciding with the edges of the standard permutohedron $\Pi_4$,
 and  diagonal edges equal to the parallel translates of the segments $r_i$. The labels of  the endpoints of a diagonal edge differ by a cyclic permutation $cycl$.

 For instance, the vertices  $(2,1,4,3)$ and $(1,4,3,2)$ share a diagonal edge which is a parallel translate of $r_2=[0,(-1,3,-1,-1)]$.
  \item The $2$-faces fall into four categories (see Figure \ref{FigCycloper5labels}):
  \begin{enumerate}
    \item Virtual hexagons of type $(q_{ij}-r_i-r_j)$. These are combinatorial permutohedra $\Pi_3$.
     \item Convex quadrilaterals, faces of the standard permutohedron $\Pi_4$.

    \item Quadrilaterals of type $q_{ij}\times (-r_k)$, and
    \item Convex hexagons,  faces of the standard permutohedron $\Pi_4$.
  \end{enumerate}
\end{enumerate}
\end{ex}

\newpage

\section{Cyclopermutohedron is the universal polytope for the moduli spaces of planar polygonal linkages}\label{SecModuliPolygon}

A \textit{polygonal  $n$-linkage} is a sequence of positive numbers
$L=(l_1,\dots ,l_n)$. It should be interpreted as a collection of
rigid bars of lengths $l_i$ joined consecutively in a chain by
revolving joints.

We assume that $L$ satisfies the triangle inequality.

 \textit{A configuration} of $L$ in the
Euclidean plane $ \mathbb{R}^2$  is a sequence of points
$R=(p_1,\dots,p_{n}), \ p_i \in \mathbb{R}^2$ with
$l_i=|p_i,p_{i+1}|$, and $l_n=|p_n,p_{1}|$.

The set $M(L)$ of all  configurations modulo orientation preserving
isometries of $\mathbb{R}^2$ is \textit{the moduli space, or the
configuration space }of the  linkage $L$.

We assume
that no configuration of $L$ fits a straight line. This assumption
implies that the moduli space $M(L)$ is a closed $(n-3)$-dimensional
manifold.

\begin{dfn}\label{admiss}
A partition of the set $\{l_1,\dots ,l_n\}$ is called
\textit{admissible } if the total length of any part does not exceed
the total length of the rest.

In the terminology of paper \cite{faS}, all parts of an admissible
partition are \textit{short sequences}.

Instead of partitions of $\{l_1,\dots ,l_n\}$ we shall speak of
partitions of the set $[n]$, keeping in mind the lengths
$l_i$.
\end{dfn}

It is proven (see \cite{pan2}) that the moduli space $M(L)$  admits a structure of a regular
cell complex (denoted by $CWM^*(L)$) whose description reads as follows.

 \begin{enumerate}
   \item  The $k$-cell of the complex are labeled by (all possible) admissible cyclically ordered partition of $[n]$ into $(n-k)$ non-empty subsets. Given a cell $C$, its label is denoted by $\lambda (C)$.
   \item A closed cell $C$  belongs to the boundary of another closed cell $C'$ whenever the label $\lambda (C')$ is finer than the label $\lambda (C)$.
 \end{enumerate}

In other words, the moduli space $M(L)$ is patched of the admissibly labeled cells of the complex $CP_{n}$.

An immediate  consequence of the construction is the following theorem:

\begin{thm}\label{TheoremUniversalModuli}
For any $n$-linkage $L$, the cell complex $CWM^*(L)$ embeds as a
subcomplex in the face lattice of the cyclopermutohedron
$\mathcal{CP}_n$.\qed
\end{thm}

\section{Proofs}\label{SectProof}

Definition \ref{DEFNFacesvirtual1} says that given a virtual
polytope $K$, its faces are associated to  non-zero vectors $\xi$.
 Thus, to list all the faces of the cyclopermutohedron,  we
need to compute the face $ \mathcal{CP}_{n+1}^\xi$ for each  $\xi$.

So we fix a non-zero vector $\xi=(x_1,x_2,...,x_n)$, assuming  by
symmetry that  $$x_1 \geq x_2 \geq ...\geq x_n.$$ Next, we  assume
that the values $x_i$ appear in $m$ \textit{clusters} of equal
entries, that is,
$$x_1=x_2=...=x_{i_1}>x_{i_1+1}=x_{i_1+2}=...=x_{i_2}>x_{i_2+1}=x_{i_2+2}=... \hbox{  etc.}$$

Denote the clusters by $U_p$ , and denote their lengths by
$d_p=i_{p+1}-i_p$:
$$U_1=\{1,2,...,i_1\}, \ U_2=\{i_1+1,i_1+2,...,i_2\}, \hbox{  etc.}$$

We say that a vector $\xi$ is \textit{diagonal-free} if $\xi$ is orthogonal to none of $r_i$.
Otherwise we call $\xi$ a \textit{diagonal} vector.

We first make two elementary observations:
\begin{lemma}\label{LemmaCluster}
 Assume that $\xi$ is a diagonal vector, which means that $\xi$ is orthogonal to some $r_j$.
 Then $\xi$ is orthogonal to $r_i$ if and only if the indices $j$ and $i$  belong to one and the same cluster.\qed
\end{lemma}

\begin{lemma}\label{LemmaPerturb}
 If $\eta$ is a small perturbation of $\xi$, then the cluster scheme related to $\eta$ refines the cluster scheme related to $\xi$.
 Moreover, given  $\xi$, a suitable choice of a small perturbation $\eta$, one gets any prescribed refinement.  \qed
\end{lemma}
\bigskip

  Since for virtual polytopes,   "face of
the sum equals sum of the faces" (see Proposition \ref{ThmFacesvirtual}), we have:
$$ \mathcal{CP}_{n+1}^\xi= \left(  \sum_{i< j} q_{ij} +S - \sum_{i=1}^n  r_i\right)^\xi  =  \sum_{i< j} q^\xi_{ij}+ S + \sum_{i=1}^n (- r_i)^\xi.$$
So we first compute the faces of the summands. The latter are
either points, or (either convex or inverse to convex) line
segments.

\subsection*{Faces of the summands}
\begin{lemma}\label{LemmaQiFACE} The face of the line segment
$q_{ij}$ is either its vertex $e_i$, or the entire segment:
 $$q_{ij}^\xi=\left\{
              \begin{array}{ll}
              e_i, & \hbox{ if $i$ and $j$ belong to different clusters }  ; \\
                q_{ij}, & \hbox{ if $i$ and $j$ belong to one and the same cluster.\qed}

              \end{array}
            \right.$$

\end{lemma}

\begin{lemma}\label{LemmaNonDiadFace}  If $\xi$ is diagonal-free,   we have:
\begin{enumerate}

      \item

      Let $k(\xi)=k \in \mathbb{N}$ be such that $$x_k>\frac{x_1+...+x_n}{n}>x_{k+1}.$$
  Then the face of the line segment $r_i$ is the  point:

  $$r_i^\xi=\left\{
              \begin{array}{ll}
               R_i, & \hbox{ if index }  i\hbox{ preceeds }U_{k}, \hbox{ that is, }i<i_{k-1}+1; \\
               0,   & \hbox{  otherwise.}
              \end{array}
            \right.
  $$

  \item By Proposition \ref{ThmFacesvirtual}, for the Minkowski inverse $-r_i$, we have
  $$(-r_i)^\xi=\left\{
              \begin{array}{ll}
                -R_i, & \hbox{ if index }  i\hbox{ preceeds }U_{k}, \hbox{ that is, }i<i_{k-1}+1; \\
               0, & \hbox{ otherwise.\qed}
              \end{array}
            \right.
  $$

\end{enumerate}
\end{lemma}
\begin{lemma}\label{LemmaDiadFace}   If $\xi$ is a diagonal vector,  then
  there is a uniquely defined
      cluster $U_k$ such that $$\forall i\in U_k, \ \ \ x_i=\frac{x_1+...+x_n}{n}.$$
      Then we have:

       $$r_i^\xi=\left\{
              \begin{array}{lll}
             r_i, & \hbox{ if }  i \in U_{k}; \\
              R_i, & \hbox{ if index }  i\hbox{ preceeds }U_{k}, \hbox{ that is, }i<i_{k-1}+1; \\
                0, & \hbox{ if index }  i\hbox{ comes after }U_{k}, \hbox{ that is, }i>i_k.
              \end{array}
            \right. $$

            and therefore,
         $$(-r_i)^\xi=\left\{
              \begin{array}{lll}
             - r_i, & \hbox{ if }  i \in U_{k}; \\
              - R_i, & \hbox{ if index }  i\hbox{ preceeds }U_{k}, \hbox{ that is, }i<i_{k-1}+1; \\
                0, & \hbox{ if index }  i\hbox{ comes after }U_{k}, \hbox{ that is, }i>i_k.\qed
              \end{array}
            \right. $$

\end{lemma}

\subsection*{Proof of Proposition \ref{ProposVert}}

 Assume that $\xi=(x_1,...,x_n) \in \mathbb{R}^n$ is a diagonal-free vector  such
that $x_1>...>x_n$ (that is, all the clusters are singletons).

Assume also that the index  $k$ satisfies
$$x_k>\frac{x_1+...+x_n}{n}>x_{k+1}.$$

 Prove that the face $$\mathcal{CP}_{n+1}^\xi=
(\Pi_n -\sum_{i=1}^n r_i)^\xi=\Pi_n^\xi -\sum_{i=1}^n(r_i)^\xi$$ is
a one-point polytope that coincides with some vertex of $\Pi_n$. (n,
The face of the permutohedron is the  vertex
$$\Pi_n^\xi=(n-1, n-2,
... , 2 , 1).$$
 By Lemmata \ref{LemmaNonDiadFace} and \ref{LemmaQiFACE}, in this particular case all the faces are also points, therefore it remains to sum them up.
 We have $n-k+1$ summands:
$$ \begin{array}{cccccc}
   (n, & n-1, & n-2, & ... & 2, & 1) \\
   +(1-n, & 1 & 1 &... & 1, & 1) \\
   +(1, &1-n, & 1 &... & 1, & 1) \\
   +(1, & 1, & 1-n, & ... & 1, & 1)
 \end{array}
$$

The sum equals  $$((n-k),(n-k-1),
(n-k-2),...,3,2,1,n,(n-1),...,(n-k+1)) \in Vert(\Pi_n)$$  which is a
vertex of the standard  permutohedron.  By the symmetry property,
each of the vertices of $\Pi_n$ arises in this way.

 Note that in the framework of the proof, we dealt with the coordinates of the vertices, not with the labels.
We remind that the label is the inverse permutation. In particular,
the vertex we discussed
 is labeled by

$$(k,(k-1),
(k-2),...,3,2,1,n,(n-1),...,(k+1)).\qed $$

\bigskip

\subsection*{Proof of Theorem \ref{TheoremMain}}
Basing on Lemmata \ref{LemmaQiFACE} , \ref{LemmaNonDiadFace}, and
\ref{LemmaDiadFace}, we put labels on the faces of the
cyclopermutohedron. By symmetry property, it suffices to assume that
$$x_1 \geq x_2 \geq ...\geq x_n.$$

Before we start putting labels, observe that the face
$\mathcal{CP}_{n+1}^\xi$ depends only on the following three
conditions: (1) whether $\xi$ is diagonal or not, (2) the cluster
scheme, and (3) the number $k=k(\xi)$ defined in Lemmata
\ref{LemmaNonDiadFace} and \ref{LemmaDiadFace}.

\subsection*{Putting labels on diagonal-free faces of $\mathcal{CP}_{n+1}$}
For a diagonal-free vector $\xi$, we say that
$\mathcal{CP}_{n+1}^\xi$ is a\textit{ non-diagonal face}, and label
the face as follows:
\begin{enumerate}
    \item Start with the label $(U_1,U_2,...,U_m)$, where $U_i$ are the clusters associated with the
     vector $\xi$. We have a linearly ordered partition of $[n]$.

    \item Add the one-element set $\{n+1\}$ \textbf{right after} $U_{k-1}$, and close the string in a circle. We get the label
$$(U_1,U_2,...,U_{k-1},\{n+1\},U_k,U_{k+1},...,,U_m)=$$
$$=((U_k,U_{k+1},...,,U_m,U_1,U_2,...,U_{k-1},\{n+1\}),$$
which should be interpreted as a cyclically ordered partition of
$[n+1]$.
\end{enumerate}

\subsection*{Putting labels on diagonal faces of $\mathcal{CP}_{n+1}$}
For a diagonal vector $\xi$, we say that $\mathcal{CP}_{n+1}^\xi$ is a\textit{ diagonal face}, and label the face
as follows:
\begin{enumerate}
    \item Start with the label $(U_1,U_2,...,U_m)$, where $U_i$ are the clusters associated with the
     vector $\xi$.

    \item Add the entry $(n+1)$   \textbf{in the set} $U_{k}$, and close the string in a circle. We get the label
$$(U_1,U_2,...,U_{k-1},\{(n+1)\}\cup U_k,U_{k+1},...,,U_m),$$
which should be interpreted as a cyclically ordered partition of
$[n+1]$.
\end{enumerate}

The following lemmata follow more or less directly from the
construction, and complete the proof of Theorem \ref{TheoremMain}.

\begin{lemma}\label{LemmaFaceOfFace}
Let a vector $\eta$ be a small perturbation of $\xi$. Then the label
of $\mathcal{CP}_{n+1}^\eta$ refines the label of
$\mathcal{CP}_{n+1}^\xi$.\qed

\end{lemma}
\begin{lemma}\label{LemmaCoincide}Two vectors $\xi$ and $\xi '$  generate one and the same cyclic label if and only if the
faces $\mathcal{CP}_{n+1}^\xi$  and $\mathcal{CP}_{n+1}^{\xi'}$
coincide.
\end{lemma}
Proof. By the label, one easily distinguishes whether $\xi$ is
diagonal or not. As we know from above, each face equals weighted
sum of (some of) the segments $q_{ij}$ and $r_i$. Which ones should
be taken, is  understood from the cluster structure. This means that
the faces with equal labels differ by translation. To prove that
they are equal, it remains to observe that
 equal labels yield equal vertex sets.\qed

\begin{lemma}\label{LemmaFaceOfFace1}
Let $\lambda$ be the label of a face $F$, and $\lambda'$ be its
refinement. There exist a vector $\xi$ and its  small perturbation
$\eta$  such that $F=\mathcal{CP}_{n+1}^\xi$, and $\lambda'$ is the
the label of $\mathcal{CP}_{n+1}^\eta$.\qed

\end{lemma}

\newpage

\section{Virtual polytopes}\label{SecVIPs}

Virtual polytopes appeared in the literature as useful
geometrization of Minkowski differences of convex polytopes. A
detailed discussion can be found in
\cite{pkh,panina:classicalProblems:2002,pan1}.

\subsection*{Virtual polytopes: definitions}
\label{ssec:groupVP}

\textit{A convex polytope }is the convex hull of a finite, non-empty point
set in the Euclidean space $\mathbb{R}^n$. Degenerate polytopes are
also included, so a closed segment and a point are polytopes, but
not the empty set. We denote by $\mathcal{P}^+$ the set of all
convex polytopes.

\begin{dfn}
Let $K$ and $L \in \mathcal{P}^+$ be  two convex polytopes.  Their \textit{Minkowski sum} $K + L$ is defined by:
$$
K + L = \{\textbf{x} +\textbf{y} : \textbf{x} \in K, \textbf{y} \in L\}.
$$
\end{dfn}

Minkowski addition   the set $\mathcal{P}^+$ to  a commutative
semigroup whose unit element is the
convex set containing exactly one point $E= \{ 0 \}$.

\begin{dfn}
The group $\mathcal{P}$ of {\em virtual polytopes} is the Grothendick group associated to the semigroup $\mathcal{P}^+$ of convex polytopes under Minkowski addition.

The elements of $\mathcal{P}$ are called {\em virtual
polytopes}.
\end{dfn}
More instructively,  $\mathcal{P}$ can be explained as follows.

\begin{enumerate}
  \item A virtual polytope is
 a formal difference $K- L$.

  \item Two such
expressions $K_1- L_1$ and $K_2- L_2$ are
identified, whenever $K_1+ L_2=K_2+ L_1$.
  \item  The group operation is
 defined by  $$(K_1- L_1) + (K_2- L_2):=
(K_1 + K_2)- (L_1 + L_2).$$

\end{enumerate}

\begin{figure}[h]
\centering
\includegraphics[width=7 cm]{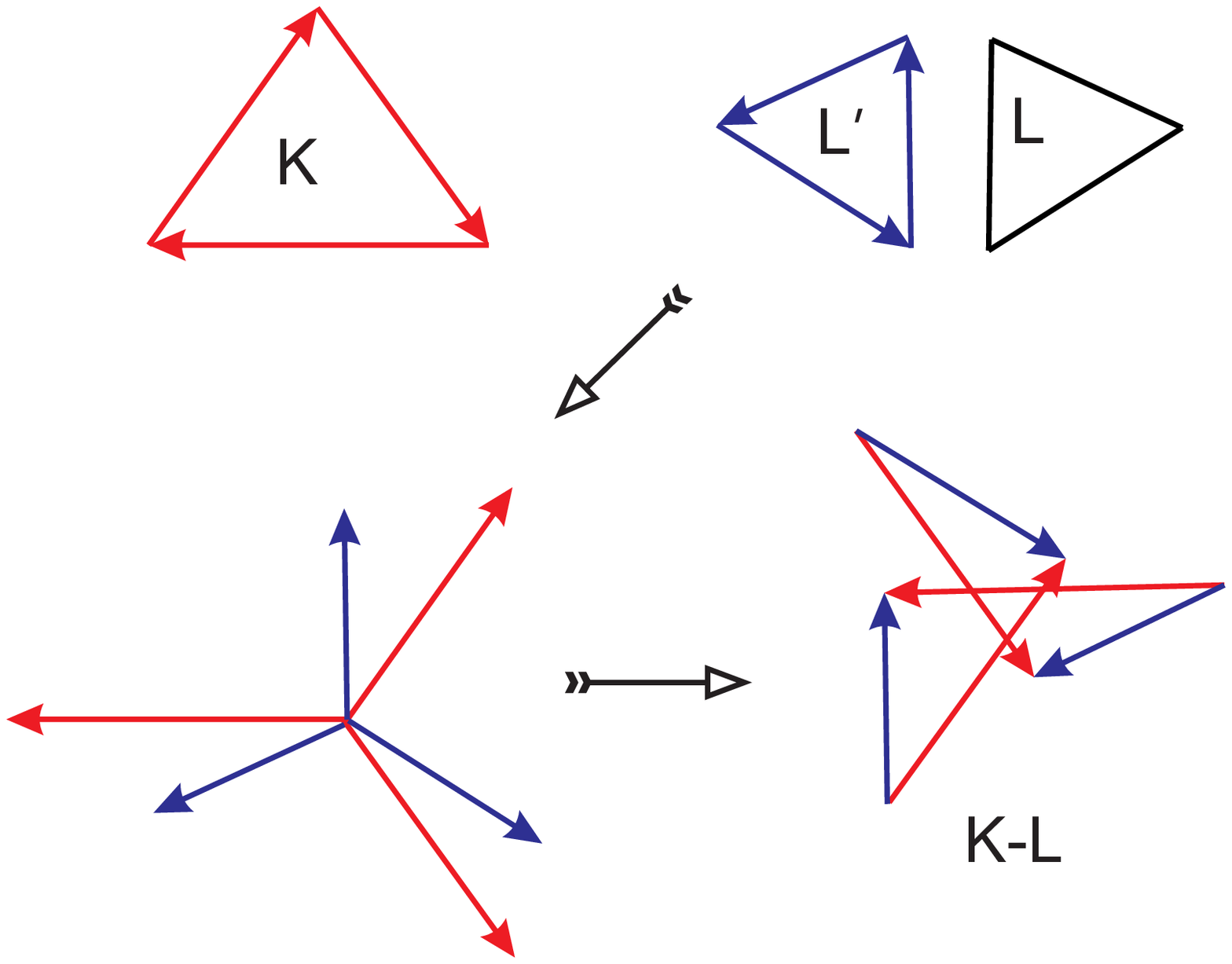}
\caption{Minkowski difference of two polygons }\label{FigMinus}
\end{figure}

\begin{figure}[h]
\centering
\includegraphics[width=7 cm]{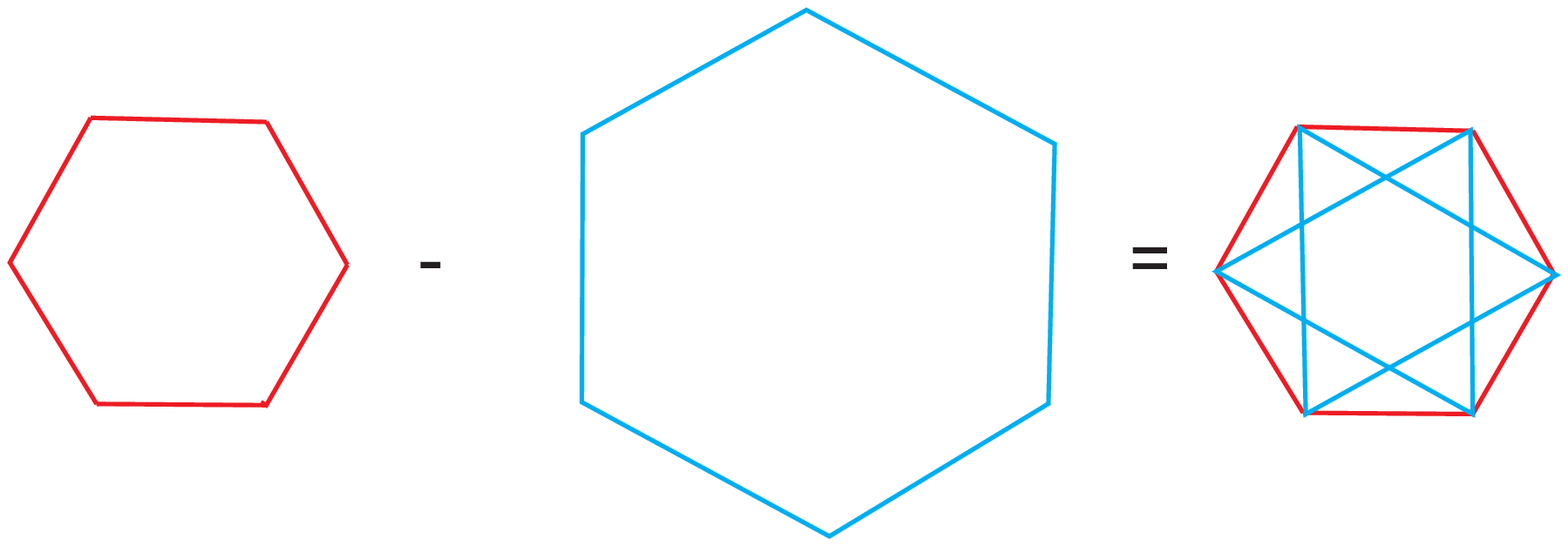}
\caption{Minkowski difference of the polygons $$  \sum_{1 \geq i> j \geq 3} q_{ij} \hbox{ and } \sum_{i=1}^3  r_i.$$  }\label{FigMinusCycloperm}
\end{figure}
\begin{alg} \label{ExMinus} In $\mathbb{R}^2$, virtual polytopes can be visualized as closed polygonal chains according a simple algorithm.
For polygons with mutually non-parallel edges (this is what we have
in the definition of cyclopermutohedron), it runs as follows.
\begin{enumerate}
    \item To compute $K-L$, orient the polygons $K$ and the symmetric image $L'$ of $L$ clockwise and counterclockwise respectively.
    \item Take the (oriented) edges of $K$ and $L'$ apart and order them    by the slope.
    \item Following the ordering, compose a closed polygonal chain, as is shown in Figure \ref{FigMinus}.
\end{enumerate}
\end{alg}

\begin{ex} \label{ExMinusCycloper}
For the particular case of Example \ref{ExCycloper4},
the algorithm gives us Figure \ref{FigMinusCycloperm}.
Observe that if we generically  perturb one of the summands, the sum would have (instead of six) twelve vertices.

\end{ex}

\subsection*{Faces of virtual polytopes}
\label{ssec:facesVirtualPolytopes}

Similar to convex polytopes, the virtual polytopes also have  well
defined faces together with incidence relations.

To explain it,we first define in an appropriate way faces of convex
polytopes.

\begin{dfn}\label{DefnFaces} Let $K \subset \mathbb{R}^n$ be a convex polytope, and
$\xi$ be a  vector in $\mathbb{R}^n$. The  polytope $K^\xi$ is the
subset of $K$ where the scalar product $(\cdot,\xi)$ attains its
maximal value:
$$K^\xi=\{x \in K | \ \forall y \in K, \ (x,\xi)\geq (y,\xi)\}.$$
 $K^\xi$  is called the \textit{face of
 $K$ with the normal vector $\xi$}, or the  \textit{face of
 $K$ in direction $\xi$}. If $\xi$ is a non-zero vector, we get a
 \textit{proper} face of the polytope.
\end{dfn}

\begin{prp}\label{ThmFacesConvex}\begin{enumerate}
                                           \item A face of a convex polytope is a convex polytope.
                                           \item A face of a face of $K$ is again a face of $K$.
                                           \item A face in the direction $\xi$ of a Minkowski sum of two convex polytopes is the Minkowski sum of the faces in direction $\xi$ from the summands:

 $(K + L)^\xi = K^\xi + L^\xi$.
                                         \end{enumerate}

\end{prp}

The item (3) of the above proposition motivates us to define faces
for virtual polytopes:

\begin{dfn}
    \label{DEFNFacesvirtual1}
Given a virtual polytope $K=L - M \in \mathcal P$, where $K$ and $L$ are
convex polytopes, and $\xi \in \mathbb{R}^n$,  the face $K^\xi$ is defined by $$K^\xi=(L -
M)^\xi :=L^\xi - M^\xi.$$
\end{dfn}

As in the convex case, we rank the faces by their dimensions.
The $k$-dimensional faces are called  \textit{ $k$-faces}; the $0$-faces  and the $1$-faces are called \textit{vertices}  and \textit{ edges} respectively.

The faces of virtual polytopes satisfy properties similar to the convex case:
\begin{prp}\label{ThmFacesvirtual}
    \cite{panina:classicalProblems:2002}\begin{enumerate}
                                          \item A face of a virtual polytope is a  virtual polytope.
                                          \item A face of a face of a virtual polytope $K$ is again a face of $K$.
                                          \item \textbf{"Face of the sum is sum of faces"}.
                                          $$(K + L)^\xi = K^\xi + L^\xi \hbox{, and }( -K)^\xi = -(K^\xi).\qed$$
                                        \end{enumerate}

\end{prp}

\begin{prp}\label{ThmFacesOf Face}
   For a virtual polytope $K$ and its face $F$, the faces of $F$ can be obtained in the following way.
   \begin{enumerate}
    \item Take all possible $\xi$ such that $F=K^\xi$.
    \item Perturb $\xi$, that is, take all vectors $\eta$ that are sufficiently close to $\xi$.
    \item For all those $\eta$, take $K^\eta$. All the resulted virtual polytopes give the set of faces of $F$.
   \end{enumerate}
\end{prp}
Proof. The statement is true for convex polytopes and extends by linearity to virtual polytopes as well.\qed

\end{document}